\documentclass[a4paper, 11pt]{article}

\usepackage[T2A]{fontenc}
\usepackage[utf8]{inputenc}
\usepackage[russian]{babel}

\usepackage{amssymb,amsmath,amsthm,comment}
\usepackage[unicode]{hyperref}
\usepackage{epstopdf}
\usepackage{graphicx}

\usepackage[left=1.5cm,right=1.5cm,top=2cm,bottom=2cm]{geometry}

\newtheorem{theorem}{Теорема}
\newtheorem*{lemma}{Лемма}
\theoremstyle{definition}
\newtheorem{definition}{Определение}
\newtheorem*{problem}{Задача}
\newtheorem*{question}{Вопрос}
\newtheorem*{conjecture}{Гипотеза}

\newcommand{\R}{\mathbb{R}}

\newcommand{\Z}{\mathbb{Z}}

\newcommand{\del}{\smash{\mskip3mu\lower1truept\hbox{\vdots}\mskip3mu}}

\renewcommand{\phi}{\varphi}

\renewcommand{\kappa}{\varkappa}

\renewcommand{\int}{\mathrm{int}}

\newcommand{\RP}{\mathbb{RP}}

\newcommand{\ver}{\vspace{.5em}}

\begin{document}

\title{Погружения многообразий}
\author{Андрей Рябичев}
\date{}

\maketitle

\begin{abstract}
Этот текст основан на прочитанном мною докладе на конференции ``Dark geometry fest'',
посвящённой геометрическим методам и их приложениям, 17 июля 2022, см.~\cite{dgf}.
Мы будем двигаться в направлении теоремы Смейла-Хирша, а для этого
разберём самый простой её частный случай --- теорему Уитни-Грауштейна о гладких кривых на плоскости.

Текст написан доступно и местами неформально, он рассчитан в первую очередь на людей,
интересующихся математикой, но пока не изучивших её достаточно глубоко.
\end{abstract}

\section{Вместо предисловия: как выучить топологию}

Основная цель данной статьи --- продемонстрировать читателю один из разделов современной математики,
называемый {\it дифференциальной топологией} (который не так широко известен, как, например, теория множеств или классический анализ).
В нём гибкие топологические методы применяются к многообразиям и отображениям между ними.
Что это такое --- об этом речь пойдёт далее.
Но чтобы не загонять себя в тупик абстрактных определений, мы начнём со вполне наглядного результата
(и попутно разберём все определения, нужные, чтобы его сформулировать и доказать).

Тем, кто сильно загорелся выучить топологию и уже жаждет выяснить путь к ней, я рекомендую для этого книги
Хатчера~\cite{hatcher} и Фоменко-Фукса~\cite{ff};
в них похожий материал изложен в немного разной последовательности и существенно разном стиле.
Чтобы углубиться в дифференциальные методы, полезен учебник Хирша~\cite{hirsch}.
Для плотного изучения буквально предмета статьи (и последующих широко обобщающих его методов)
я могу порекомендовать книги Громова~\cite{gromov-book} и Мишачёва-Элиашберга~\cite{eliashberg-mishachev-book} ---
они, опять же, написаны в совершенно разном стиле, и первая далеко покрывает содержание второй.

С моей точки зрения, однако, главное в изучении математики
--- вовсе не чтение большого количества хорошо подобранных книжек и статей,
а в первую очередь --- самостоятельное решение задач и их сдача преподавателю.
В этом желающим может помочь, например, Независимый московский университет.
И функция эта не может быть реализована нашей статьёй, цель которой ---
возвращаясь к началу, --- лишь обозначить направление исследования.

\section{Кривые на плоскости}

\begin{definition}
{\it Замкнутой кривой на плоскости} называется непрерывное отображение $\alpha:S^1\to\R^2$.
\end{definition}

Здесь $S^1$ обозначает окружность, $\R^2$ --- плоскость.
Под {\it отображением} из окружности в плоскость мы подразумеваем любое правило (соответствие),
при котором у каждой точки $x$ на окружности есть ровно один образ $\alpha(x)$ на плоскости.
Образы разных точек окружности могут совпадать на плоскости, это не запрещено.

\ver\hfil
\includegraphics[height=25mm]{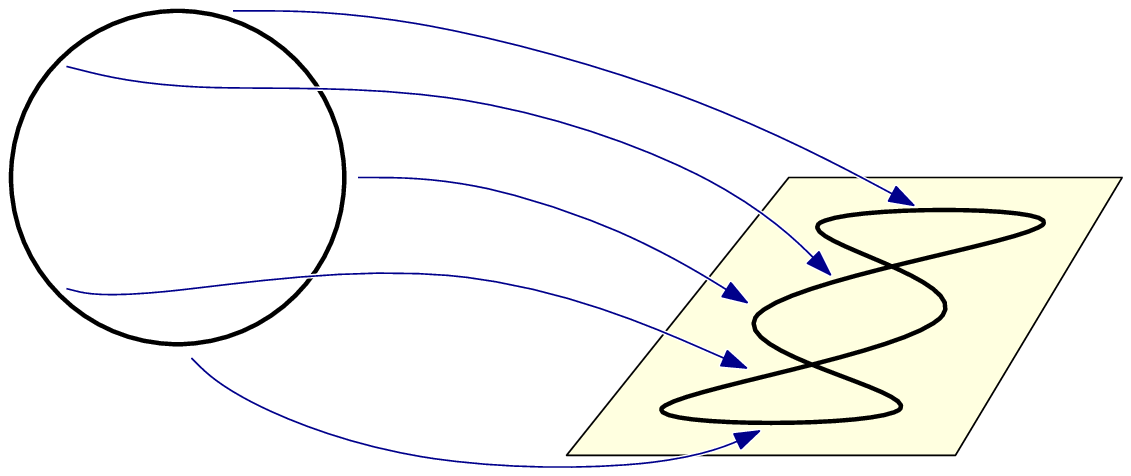} \hfil
\includegraphics[height=25mm]{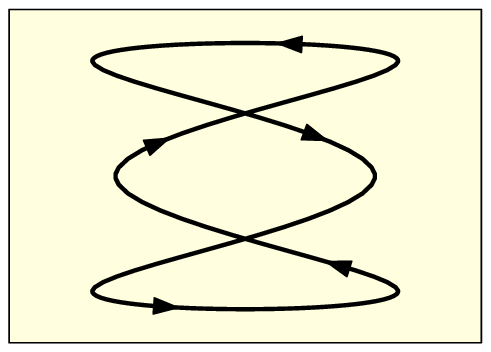}
\ver

{\it Непрерывность} нашего отображения проще воспринимать умозрительно,
так что {\it близкие точки переходят в близкие};
исчерпывающая формализация этого понятия чересчур трудоёмка и требует отдельной большой работы, выходящей далеко за рамки статьи.

Чтобы задать кривую в тексте, мы часто будем рисовать лишь множество точек на плоскости, которым что-либо соответствует (образ отображения).
Хотя, формально, соответствие между множеством точек на окружности и нарисованным множеством
несёт гораздо больше информации --- о том, какая именно точка куда перешла.
Как мы увидим впоследствии, почти всей этой информацией мы для наших целей можем пренебречь.
Почти всей --- за исключением ``направления обхода'' кривой на плоскости,
по мере того как точка на окружности бежит в положительном направлении (против часовой стрелки).

\section{Гладкие кривые}

\begin{definition}
Мы будем называть кривую $\alpha:S^1\to\R^2$ {\it гладкой},
если её вектор скорости $\dot\alpha$ в каждой точке $x\in S^1$ ненулевой
и $\dot\alpha$ непрерывно зависит от $x$.
\end{definition}

На рисунке ниже кривая, изображённая слева, является гладкой, а кривая справа --- нет,
поскольку в ``точке излома'' у неё нельзя определить вектор скорости (производную), чтобы он был ненулевым.%
\footnote{Вообще говоря, {\it гладкие кривые}, по общепринятому определению гладкости, должны были бы являться лишь {\it непрерывно-дифференцируемыми},
а кривые с ненулевым вектором скорости правильнее было бы называть {\it погруженными кривыми}.
Используемая нами терминология иногда применяется именно по отношению к кривым,
к тому же она в целом выглядит более лаконичной и менее пугающей.}

\ver\hfil
\includegraphics[height=40mm]{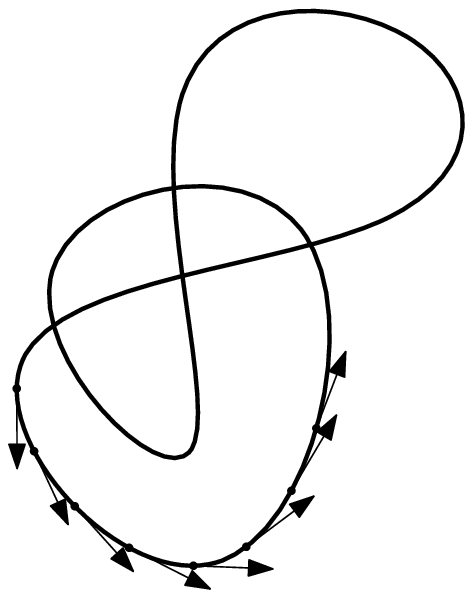} \hfil
\includegraphics[height=40mm]{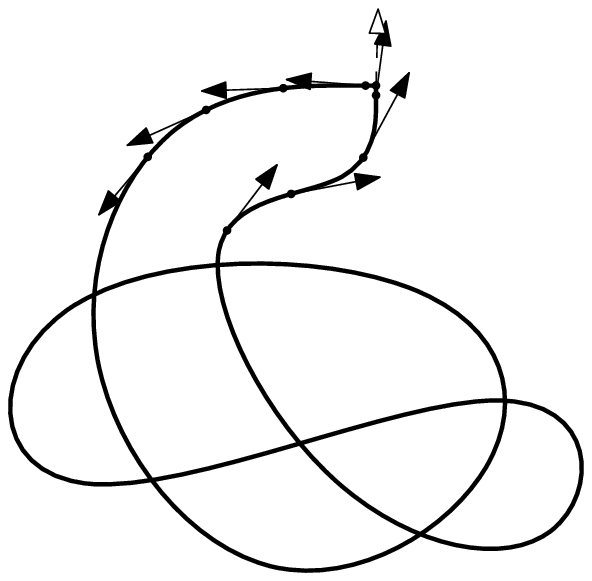}
\ver

Другими словами, $\dot\alpha$ задаёт непрерывное отображение $S^1\to\R^2\setminus(0,0)$.
Действительно, вектор скорости $\dot\alpha$ может принимать любые значения на плоскости, кроме начала координат $(0,0)$.

Определив новый объект, естественно попытаться решить задачу о классификации:

\begin{question}
Как описать все гладкие замкнутые кривые на плоскости?
\end{question}

\ver\hfil
\includegraphics[height=18mm]{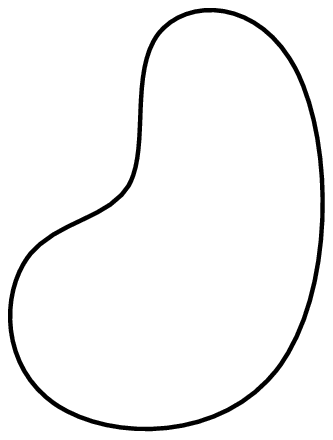} \hfil
\includegraphics[height=18mm]{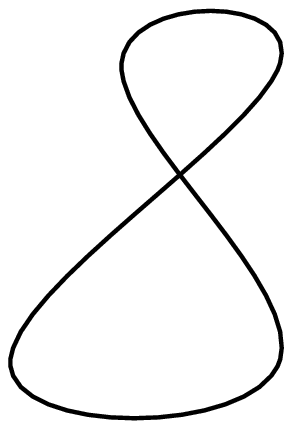} \hfil
\includegraphics[height=18mm]{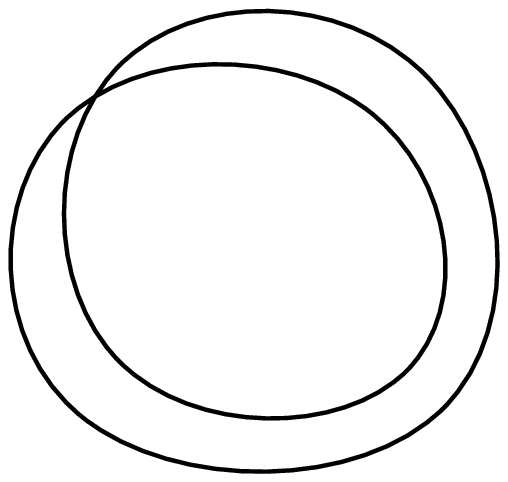} \hfil
\includegraphics[height=18mm]{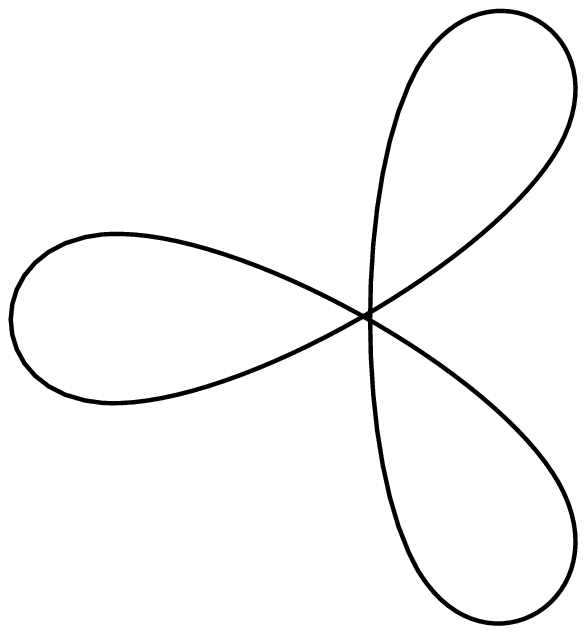} \hfil
\includegraphics[height=18mm]{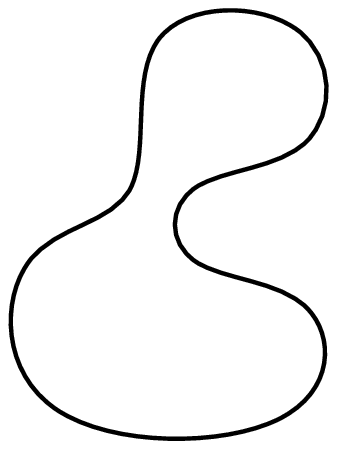} \hfil
\includegraphics[height=18mm]{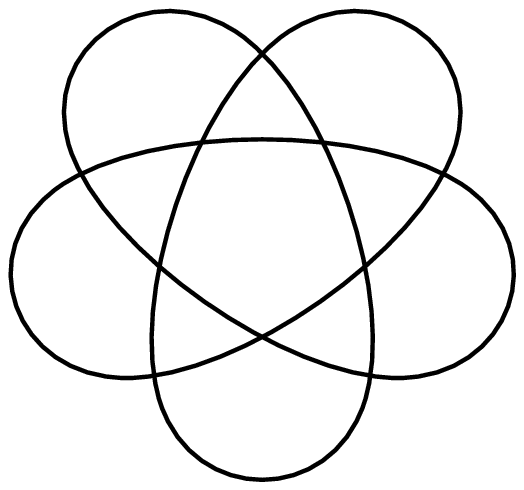} \hfil
\includegraphics[height=18mm]{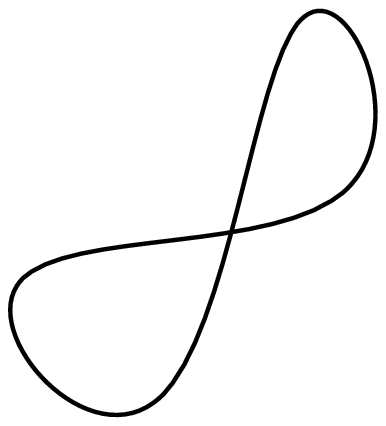}
\ver

Понятно, что всего таких кривых бесконечно много --- более точно, континуальное множество.
Было бы удобнее решать задачу о классификации, считая некоторые кривые ``эквивалентными'' --- но какие именно?

Например, пятая кривая на рисунке выше --- это всего лишь немного сплющенная первая.
А вторая похожа на повёрнутую и растянутую последнюю.
Это наводит на мысль о некоторых деформациях, как, собственно, и принято в топологии\ldots

\section{Деформация гладких кривых}\label{s:deform}

Гладкую кривую можно деформировать так, чтобы она оставалась гладкой в любой момент времени.

\ver\hfil
\includegraphics[height=20mm]{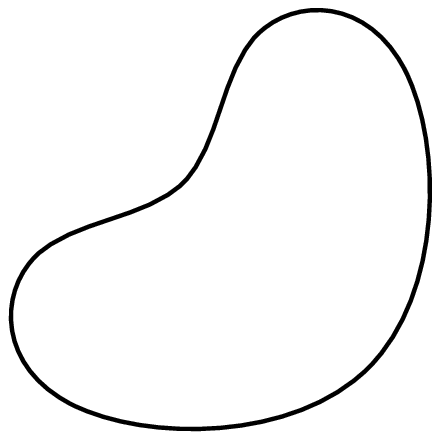} \hfil
\includegraphics[height=20mm]{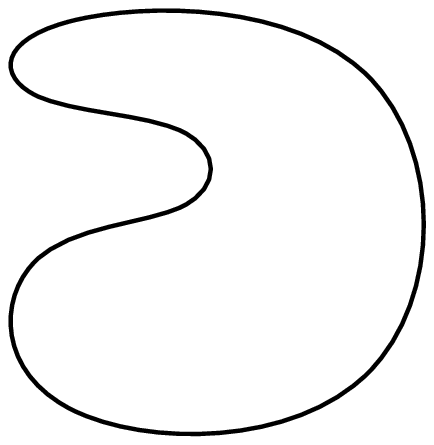} \hfil
\includegraphics[height=20mm]{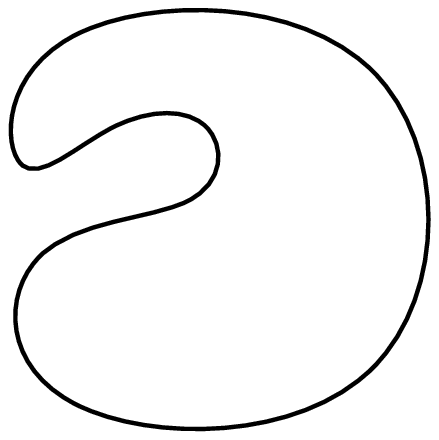} \hfil
\includegraphics[height=20mm]{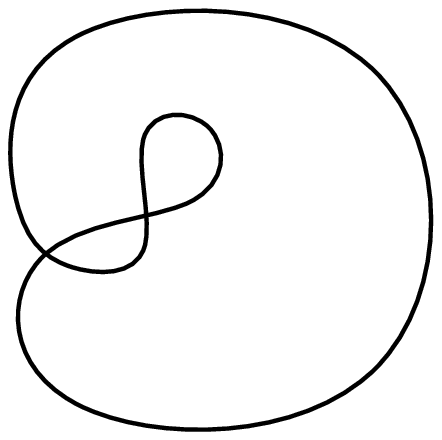} \hfil
\includegraphics[height=20mm]{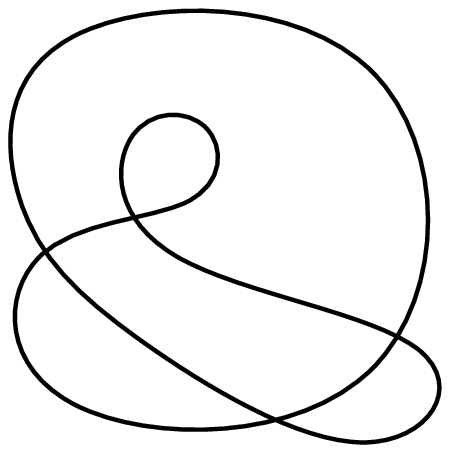}
\ver

\noindent
От рассматриваемых деформаций мы дополнительно требуем, чтобы вектор скорости в каждой точке зависел от времени непрерывно
(поскольку это тоже важная часть структуры гладкой кривой, а все преобразования должны быть непрерывными).

\ver\hfil
\includegraphics[height=11mm]{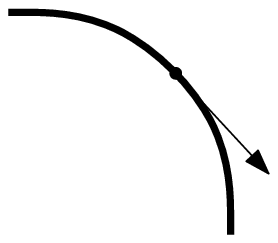} \hfil
\includegraphics[height=11mm]{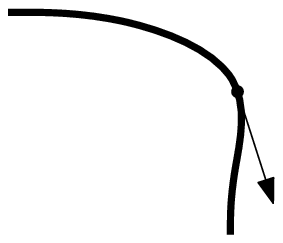} \hfil
\includegraphics[height=11mm]{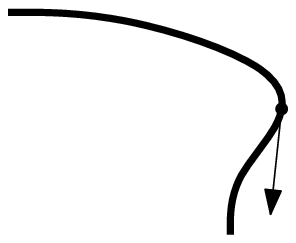} \hfil
\ver

\noindent
Деформация, удовлетворяющая этим двум условиям, называется {\it регулярной гомотопией}.
Про кривые, которые можно получить друг из друга с её помощью, говорят, что они {\it регулярно гомотопны}.

Таким образом, наше определение замкнуто на сущности гладкой кривой и не использует никаких сторонних понятий, поэтому оно довольно естественно.
Важно отметить при этом, что в процессе регулярной гомотопии
кривая может не только растягиваться и сжиматься, но ещё на ней могут появляться и исчезать самопересечения (как в примере выше).

Более формально, регулярная гомотопия это гладкое отображение $S^1\times[0;1]\to\R^2$.
Следующее утверждение несложно доказать, пользуясь любой из версий определения регулярной гомотопии.

\begin{lemma}
Регулярная гомотопность является отношением эквивалентности.
\end{lemma}

Таким образом, все кривые делятся на непересекающиеся {\it классы эквивалентности},
которые в дальнейшем и будут представлять для нас интерес.

\section{Примеры регулярных гомотопий}

Возьмём такие две кривые $\alpha_0$ и $\alpha_1$:

\ver\hfil
\includegraphics[height=35mm]{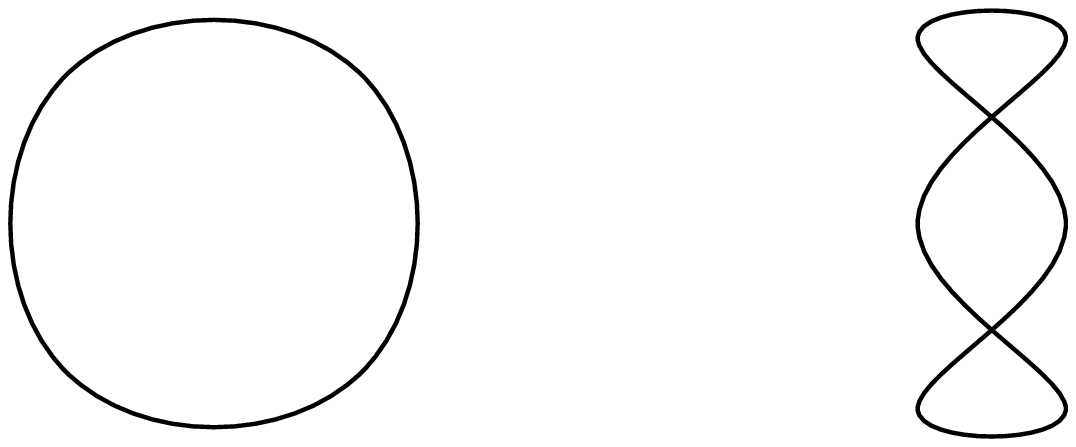}
\ver

Существует ли между ними регулярная гомотопия?
Немного подумав, мы получаем утвердительный ответ, построив регулярную гомотопию явно:

\ver\hfil
\includegraphics[height=20mm]{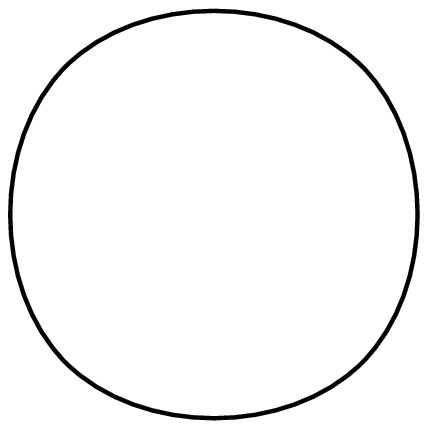} \hfil
\includegraphics[height=20mm]{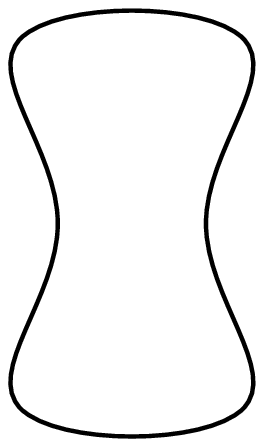} \hfil
\includegraphics[height=20mm]{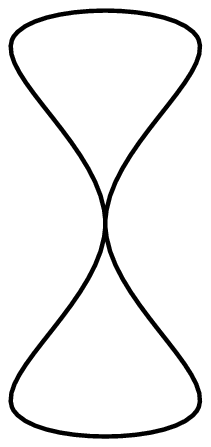} \hfil
\includegraphics[height=20mm]{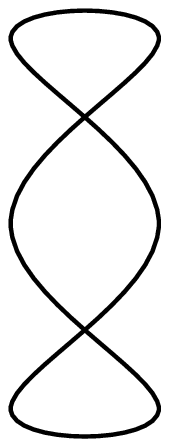}
\ver

А для таких двух кривых?

\ver\hfil
\includegraphics[height=35mm]{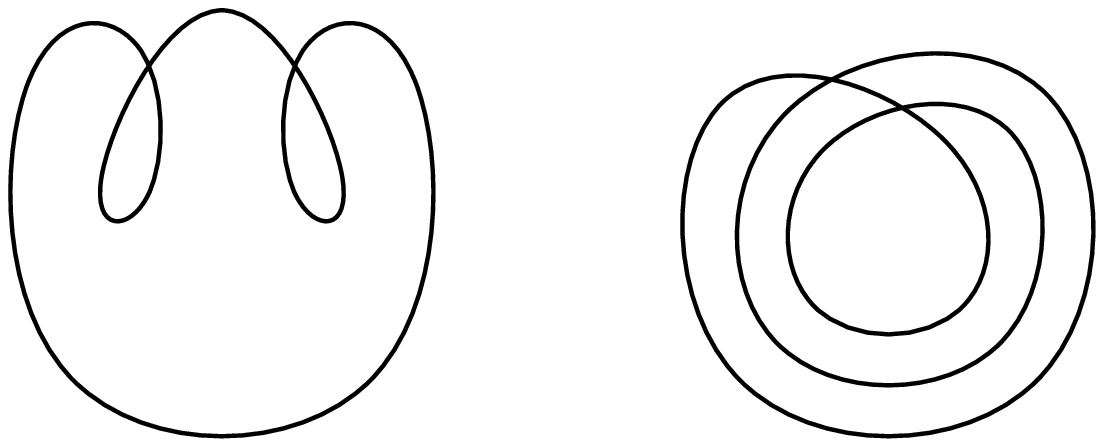}
\ver

Здесь ситуация уже кажется запутанной.
Тем не менее, мы всё же можем построить желаемую регулярную гомотопию и это не так сложно:

\ver\hfil
\includegraphics[height=20mm]{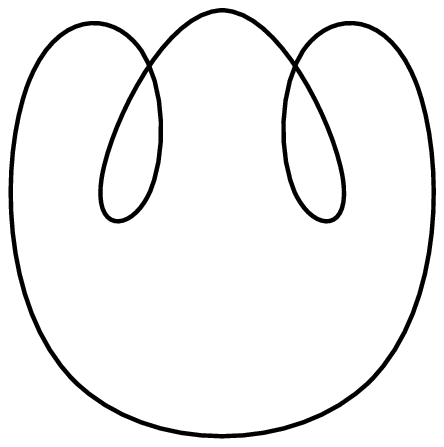} \hfil
\includegraphics[height=20mm]{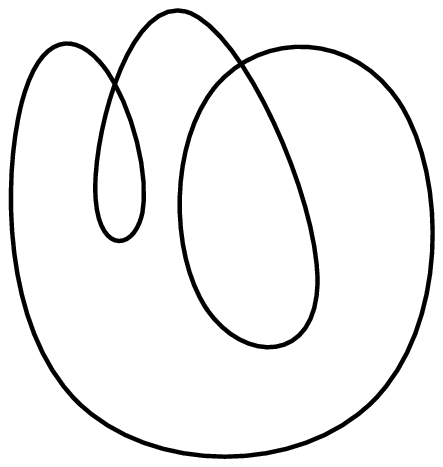} \hfil
\includegraphics[height=20mm]{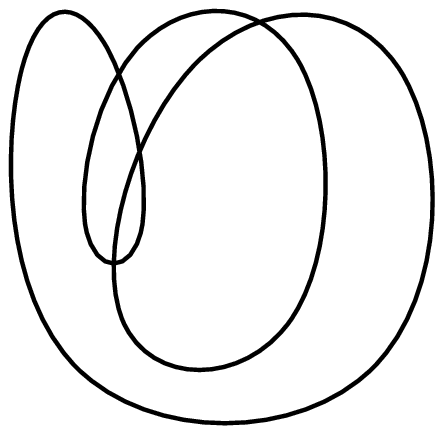} \hfil
\includegraphics[height=20mm]{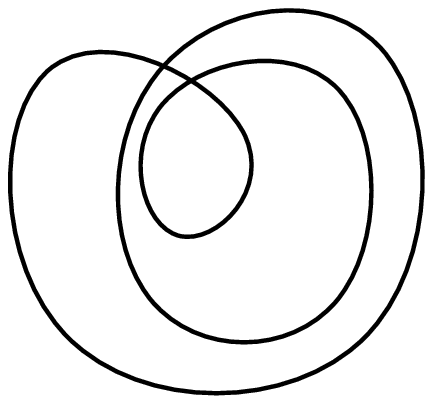} \hfil
\includegraphics[height=20mm]{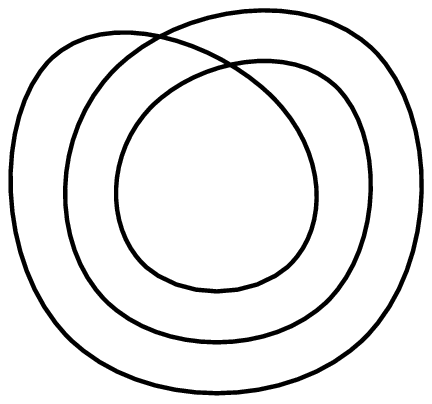}
\ver

А для таких?

\ver\hfil
\includegraphics[height=35mm]{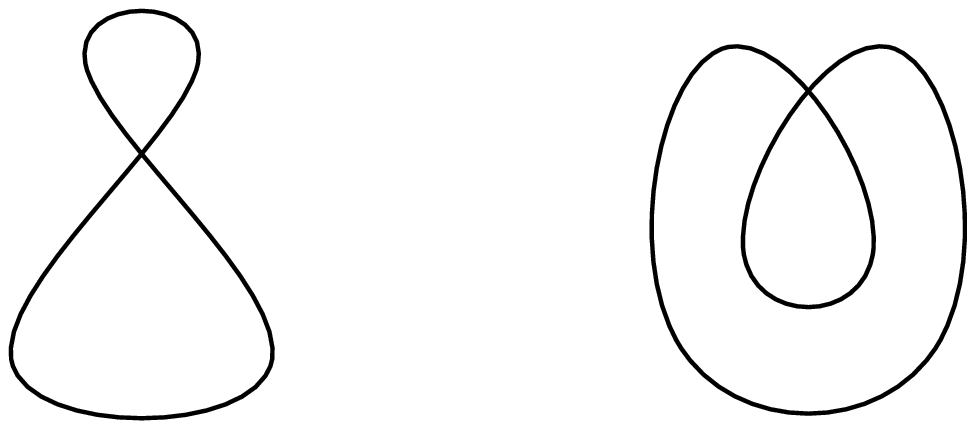}
\ver

Попытки построить регулярную гомотопию так же легко, как в предыдущие разы, встречают некое препятствие:
в какой-то момент в процессе деформации кривая перестаёт быть гладкой.
Кажется, что ответ --- ``нет'',  но как это доказать?
Ведь то, что у нас не получается предъявить регулярную гомотопию, вовсе не означает, что её не существует ---
она может быть устроена довольно сложно и т.\,п.

Чтобы решать задачи подобного сорта, в топологии (как и во многих других областях) используется один и тот же важный математический приём.

\section{Общеобразовательное отступление. Инварианты}

Чтобы доказать, что что-то нельзя сделать, часто используют {\it инварианты}.

Инвариант --- некая величина, не меняющаяся в рамках данной задачи.
Зачастую нахождение правильного инварианта равносильно решению задачи.
Для тех читателей, кто мало сталкивался с инвариантами ранее, мы разберём несколько почти тривиальных примеров.

\begin{problem}
В ряд стоят 100 фишек.
Разрешено менять местами фишки, стоящие через одну.
Можно ли переставить фишки в обратном порядке?
\end{problem}

{\bf Ответ:} {\it нет, нельзя.}
Первая фишка может оказаться только на нечётных позициях.
Однако в конце ей надо оказаться на сотой позиции, т.\,е.\ на чётной.
{\it Чётность позиции первой фишки --- инвариант.}

\begin{problem}
Можно ли квадрат со стороной 1
разрезать на части и сложить из них два правильных треугольника со стороной 1?

\ver\hfil
\includegraphics[height=12mm]{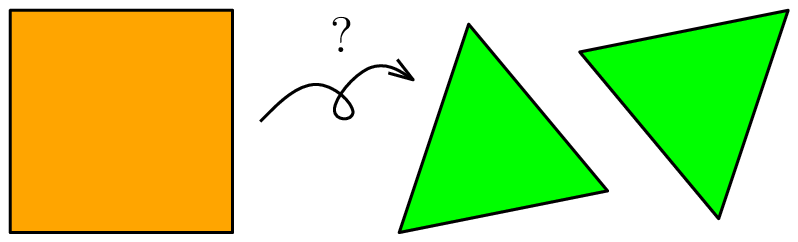}
\end{problem}

{\bf Ответ:} {\it нет, нельзя.}
Площадь квадрата больше суммарной площади двух треугольников (докажите это самостоятельно!).
{\it Площадь --- инвариант.}

\begin{problem}
Маша разломила шоколадку на 8 кусков.
Затем она разломила один из кусков еще на 8, и т.\,д.
Могло ли у неё в какой-то момент получиться ровно 179 кусков?
\end{problem}

{\bf Ответ:} {\it нет, не могло.}
Каждый раз число кусков даёт остаток 1 при делении на 7, в то время как 179 даёт остаток 4.
{\it Остаток $\pmod 7$ --- инвариант.}

\section{Инварианты регулярной гомотопии гладких кривых}

Какой же инвариант ``различает'' эти две кривые?

\ver\hfil
\includegraphics[height=25mm]{regular-homotopic-3.eps}
\ver

{\bf Ответ:} {\it число оборотов вектора скорости.}
Назовём это число $r(\alpha)$.
Чтобы его вычислить, нужно, двигаясь по кривой, отмечать направление вектора скорости,
а затем сосчитать, какое суммарное число оборотов (с учётом знака) у нас получилось.
Обороты против часовой стрелки считаются со знаком ``$+$'',
а по часовой --- со знаком ``$-$''.

Действительно, у первой кривой это число равно 0, а у второй оно равно 2.
Это можно посчитать вручную, внимательно проследив за вектором скорости при обходе кривой.

\ver\hfil
\includegraphics[height=30mm]{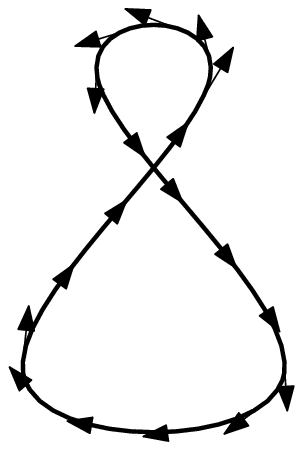}\hfil
\includegraphics[height=30mm]{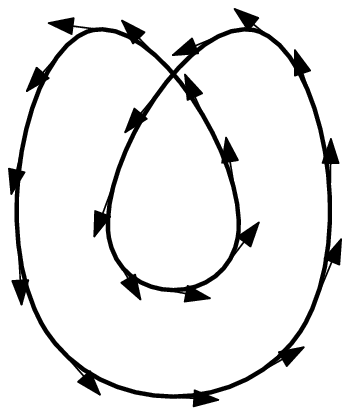}
\ver

Но как же доказать, что функция $r(\,\cdot\,)$ является инвариантом?
Для этого используется следующий аргумент:
в процессе регулярной гомотопии кривой $\alpha$
число $r(\alpha)$ тоже должно {\it меняться непрерывно}.

\ver\hfil
\includegraphics[height=30mm]{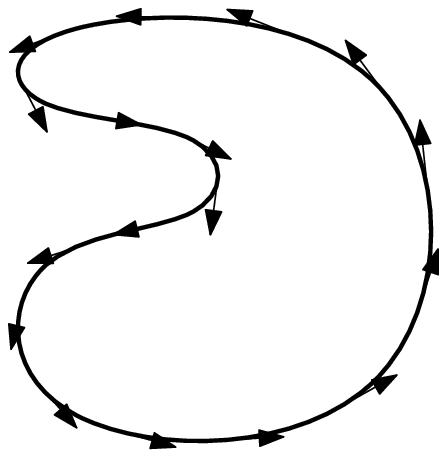} \hfil
\includegraphics[height=30mm]{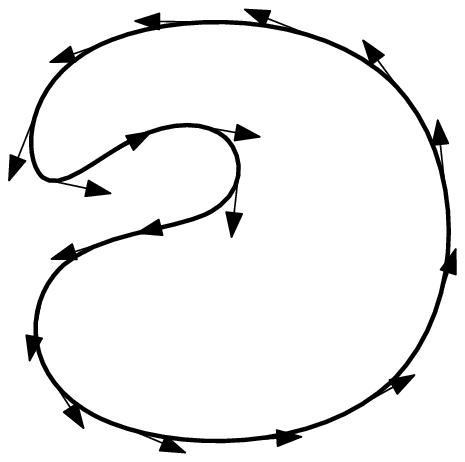} \hfil
\includegraphics[height=30mm]{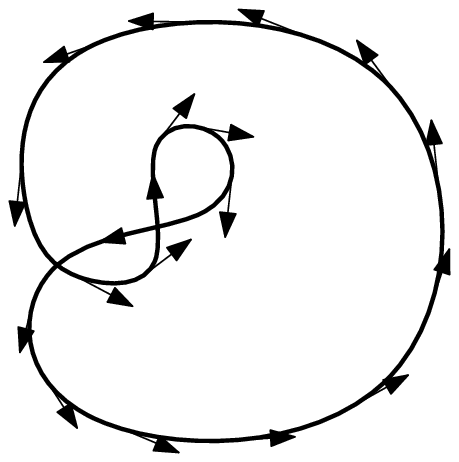} \hfil
\ver

\noindent
Но в каждый момент времени $r(\alpha)$ является целым.
Поэтому оно вообще не может меняться!
Строгое доказательство можно найти, например, в \cite{hatcher}, \S1.1.

\section{Попарно не регулярно гомотопные кривые}\label{s:collection}

Для кривой $\alpha$ число $r(\alpha)$ может принимает целое значение.
Легко придумать серию кривых, значения функции $r(\,\cdot\,)$ на которых --- в точности все целые числа.
Нужно просто вращать вектор скорости разное количество раз, например так:

\ver\hfil
\includegraphics[height=20mm]{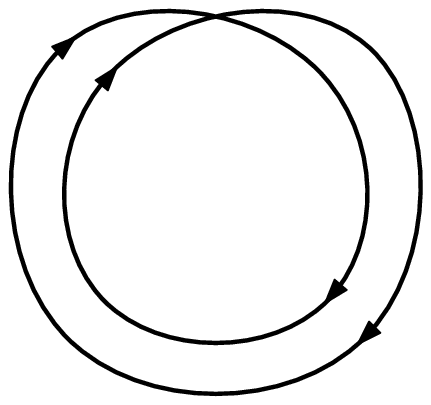}\hfil
\includegraphics[height=20mm]{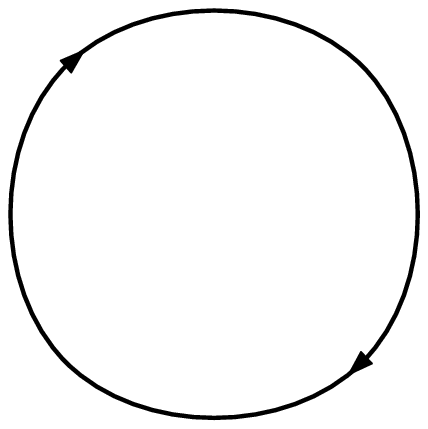}\hfil
\includegraphics[height=20mm]{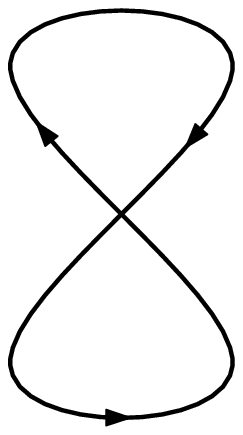}\hfil
\includegraphics[height=20mm]{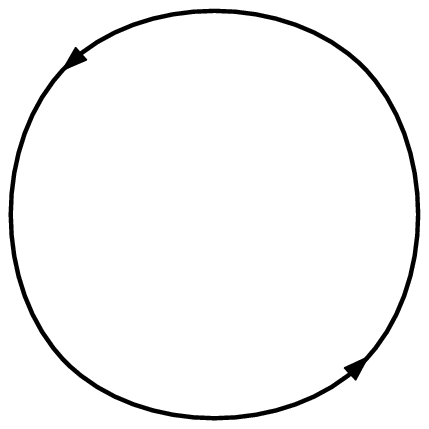}\hfil
\includegraphics[height=20mm]{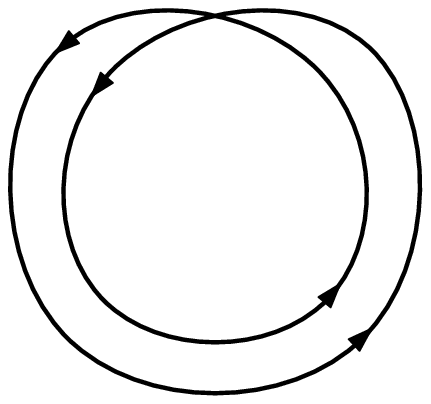}\hfil
\includegraphics[height=20mm]{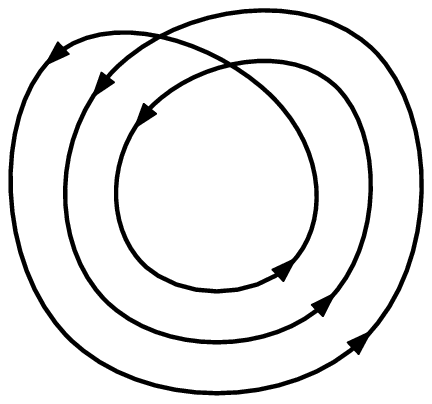}
\ver

Таким образом, мы получили бесконечную серию попарно не регулярно гомотопных кривых.
Назовём их соответственно $\ldots,\ \gamma_{-2},\ \gamma_{-1},\ \gamma_{0},\ \gamma_{1},\ \gamma_{2},\ \gamma_{3},\ldots$,
по числу оборотов вектора скорости, т.\,е. так что $r(\gamma_m)=m$.

Кстати, если пройти кривую $\gamma_0$ в противоположную сторону,
получится кривая $\gamma_0'$, вообще говоря, не входящая в наш набор.

\ver\hfil
\includegraphics[height=25mm]{curve-etalon-0.eps}\hfil
\includegraphics[height=25mm]{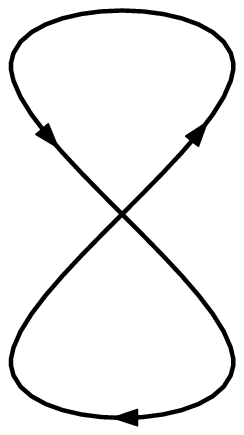}

\begin{problem}
Кривые $\gamma_0$ и $\gamma'_0$ на рисунке выше регулярно гомотопны.
\end{problem}

Далее для краткости мы будем считать регулярно гомотопные кривые {\it одинаковыми} и обращаться к ним как к одному и тому же объекту.
С этой точки зрения, для кривой $\gamma_0$ не существенно направление обхода ---
в то время как для любой другой $\gamma_m$, при $m\ne0$, изменение направления обхода превращает её в {\it другую} кривую $\gamma_{-m}$.

\section{Классификация кривых с точностью до регулярной гомотопии}

Теперь естественно задаться вопросом:
существуют ли кривые (с точностью до регулярной гомотопии), не входящие в серию
$\ldots,\ \gamma_{-2},\ \gamma_{-1},\ \gamma_{0},\ \gamma_{1},\ \gamma_{2},\ \gamma_{3},\ldots$?
Других кривых бы не существовало в случае выполнения следующего предположения:

\begin{conjecture}
Кривые с одинаковым числом $r(\,\cdot\,)$ регулярно гомотопны.
\end{conjecture}

Такое предположение согласуется с предыдущим упражнением про кривые $\gamma_0$ и $\gamma'_0$.
Действительно, мы видим, что $r(\gamma_0)=r(\gamma_0')=0$, а ещё $\gamma_0$ и $\gamma'_0$ регулярно гомотопны.

Проверим нашу гипотезу для следующей кривой $\alpha$

\ver\hfil
\includegraphics[width=35mm]{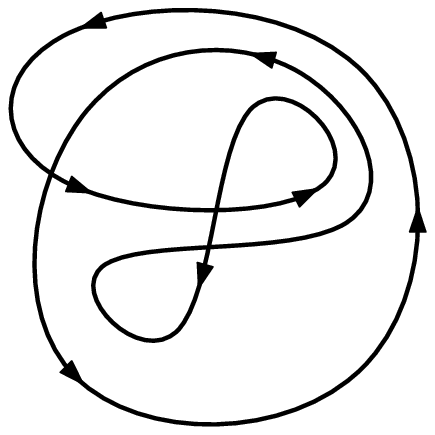}
\ver

Для неё число $r(\alpha)$ равно 2.
Значит, она должна быть
регулярно гомотопна кривой $\gamma_2$ из перечня выше.
И действительно, регулярную гомотопию можно построить явно:

\ver\hfil
\includegraphics[width=20mm]{homotopy-1.eps}\hfil
\includegraphics[width=20mm]{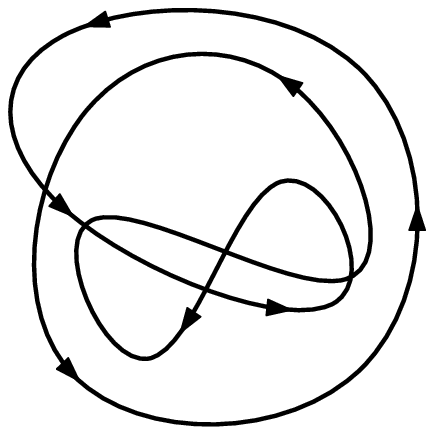}\hfil
\includegraphics[width=20mm]{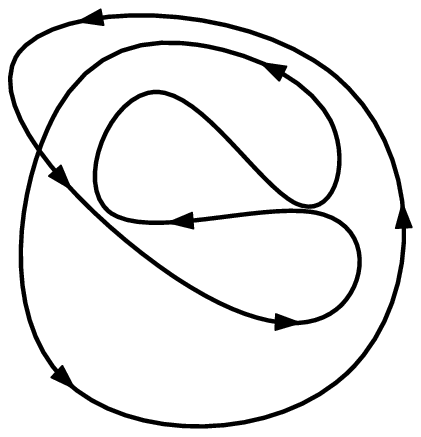}\hfil
\includegraphics[width=20mm]{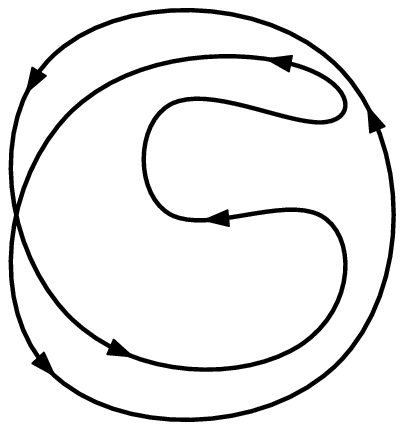}\hfil
\includegraphics[width=20mm]{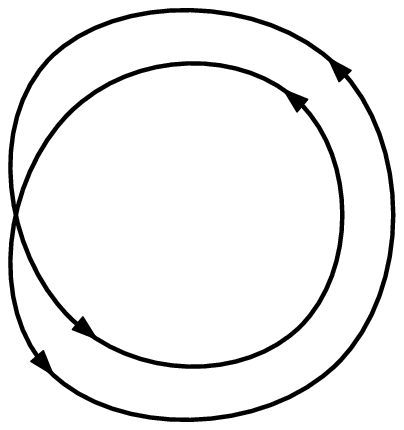}

\section{Доказательство Гипотезы о классификации}

Как можно доказать Гипотезу?
Естественная идея --- привести все кривые к неким каноническим видам!

А именно, если дана кривая $\alpha$, можно попробовать её максимально упростить,
в надежде что получится $\gamma_k$ для какого-нибудь $k$.

\subsection*{Доказательство Гипотезы. Шаг 1}

Возьмём кривую $\alpha$.

\ver\hfil
\includegraphics[width=50mm]{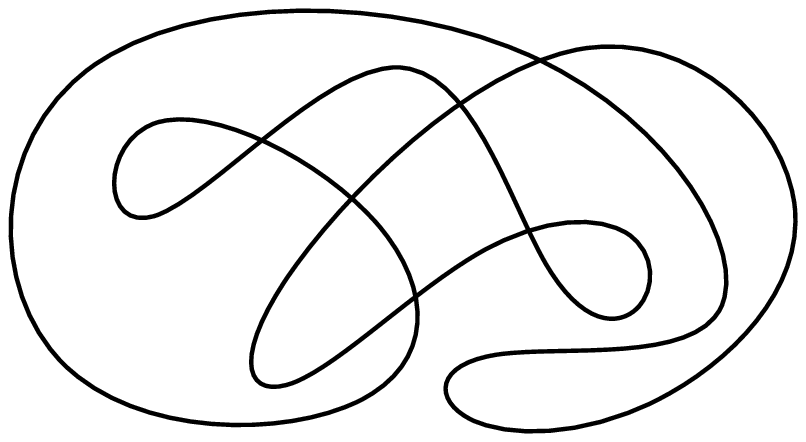}
\ver

Выберем прямую $l$, такую что $\alpha$ имеет лишь конечное число касательных, перпендикулярных $l$.
Если вектор скорости $\dot\alpha$ всё время поворачивается, то подойдёт любая прямая.%
\footnote{В общем случае, такая прямая $l$ существует по лемме Сарда.
Однако, чтобы применить лемму Сарда к $\dot\alpha$, надо потребовать для $\alpha$ гладкости порядка 2.
С другой стороны,  в геометрии часто принято считать все отображения бесконечно гладкими,
поэтому мы не должны чувствовать себя обременёнными таким дополнительным предположением.}

\ver\hfil
\includegraphics[width=60mm]{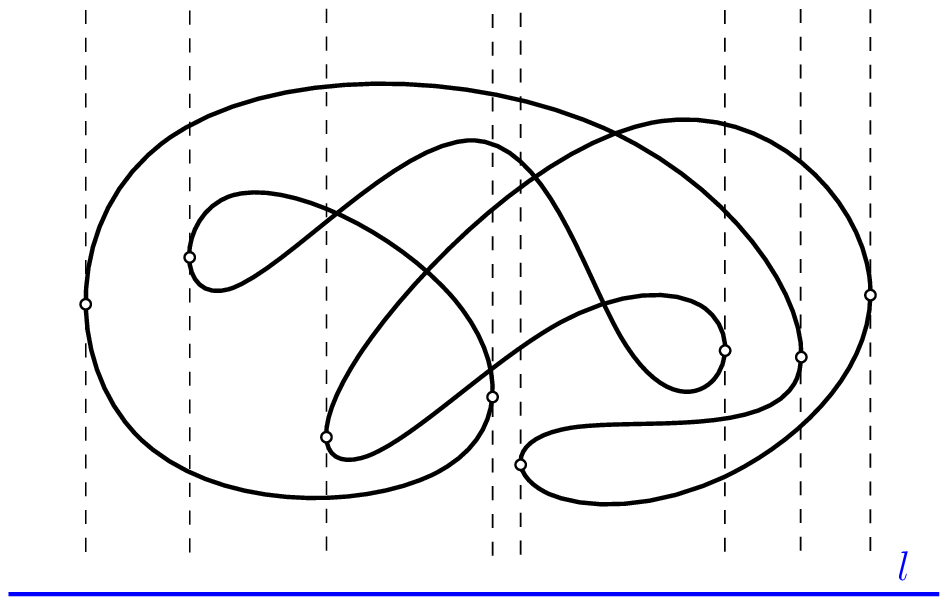}
\ver

\subsection*{Доказательство Гипотезы. Шаг 2}

Далее ``сплюснем'' кривую $\alpha$ так, чтобы она состояла
из участков, ``почти параллельных $l$'', и ``разворотов''.

\ver\hfil
\includegraphics[width=30mm]{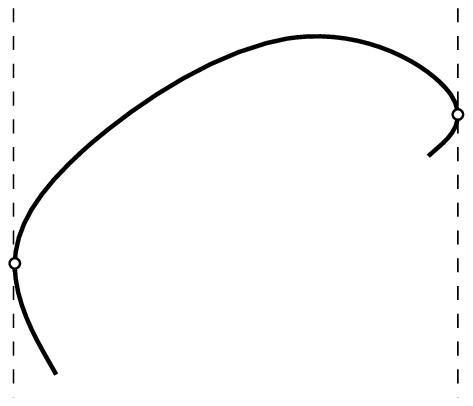}\hfil
\includegraphics[width=30mm]{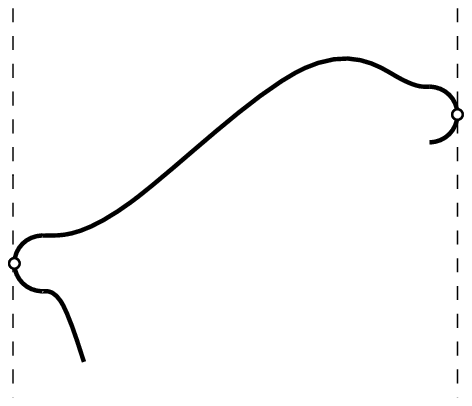}\hfil
\includegraphics[width=30mm]{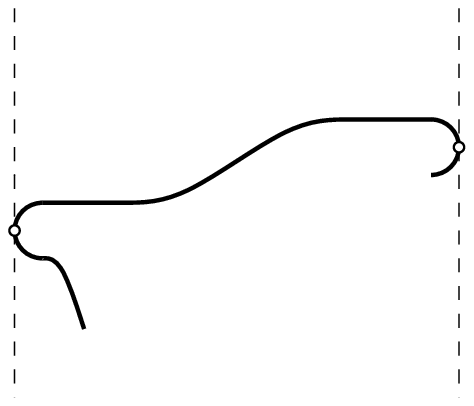}

Назовём полученную кривую $\alpha'$.
Фактически, $\alpha'$ описывается последовательностью разворотов, их направлением и их
``координатой'' по оси $l$, все остальные данные легко меняются под действием регулярной гомотопии.

\ver\hfil
\includegraphics[width=60mm]{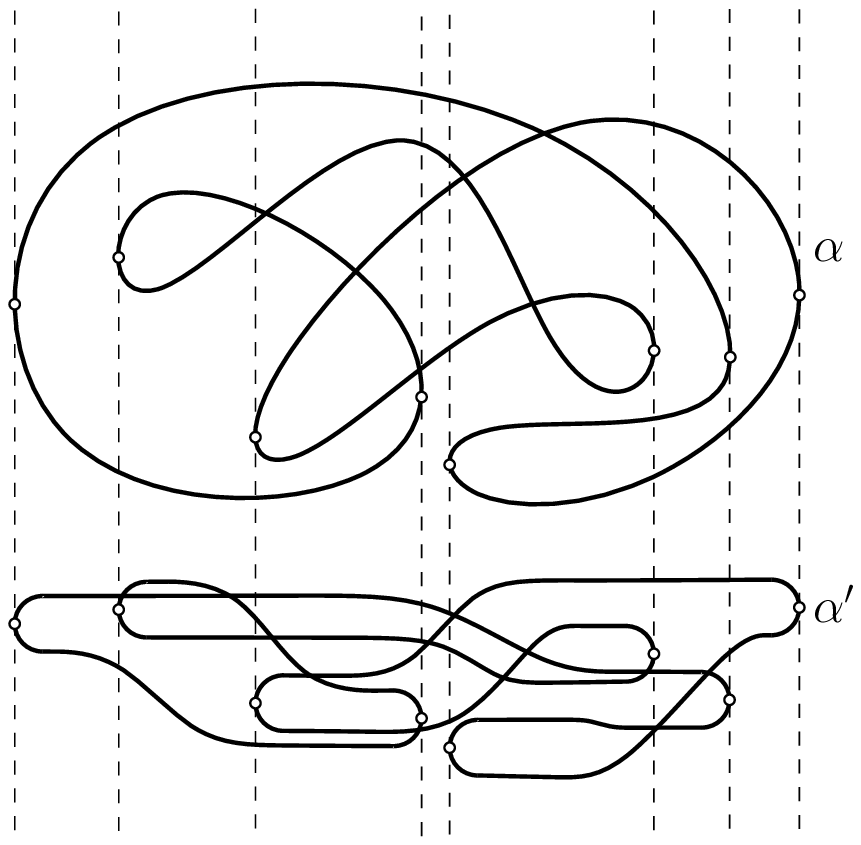}

\subsection*{Доказательство Гипотезы. Шаг 3}

Заметим, что по ходу движения можно различить развороты ``направо'' и ``налево'',
в том смысле как если бы мы сами ехали по кривой.
При этом развороты каждого из этих видов могут выглядеть двумя способами.
На рисунке ниже --- пара разворотов ``направо'' и пара разворотов ``налево''.

\ver\hfil
\includegraphics[width=27mm]{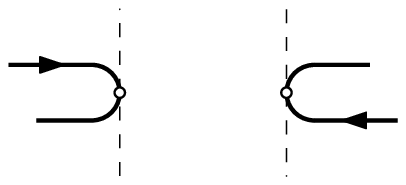}\hfil
\includegraphics[width=27mm]{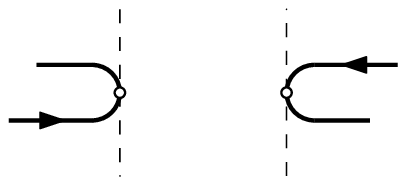}
\ver

Теперь будем деформировать $\alpha'$ следующим образом.
Если подряд встречаются два разворота в противоположных направлениях, то их можно {\it сократить}.

\ver\hfil
\includegraphics[height=13mm]{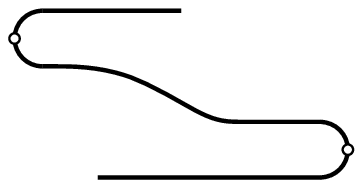} \hfil
\includegraphics[height=13mm]{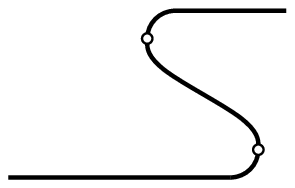} \hfil
\includegraphics[height=13mm]{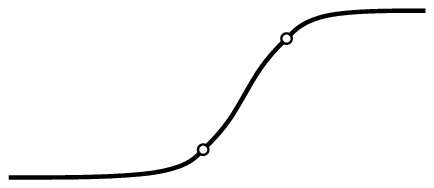}

\subsection*{Доказательство Гипотезы. Шаг 3. Уточнение}

Более точно, пусть нам встретилось два разворота в противоположных направлениях.
Перед первым разворотом и после второго есть ещё какие-то фрагменты кривой.

Если требуется, ``потянем'' за концы до первого разворота и после второго разворота
так, чтобы в итоге эти концы накладывались при проекции на $l$.
При этом на фрагментах кривой за концами все развороты остаются неизменными
и лишь ``раздвигаются'' в противоположные стороны,
новых разворотов не возникает.

\ver\hfil
\includegraphics[height=20mm]{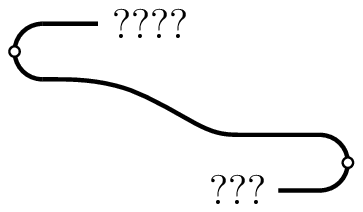}\hfil
\includegraphics[height=20mm]{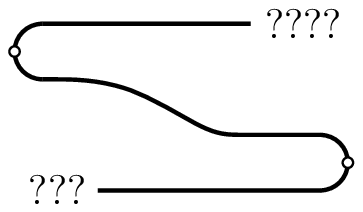}
\ver

С полученным участком кривой можно проделать процедуру {\it сокращения разворотов}
(хотя с исходным, возможно, было нельзя).

\ver\hfil
\includegraphics[height=15mm]{turns-det-killing-2.eps} \hfil
\includegraphics[height=15mm]{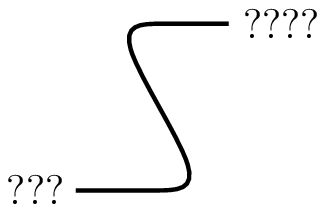} \hfil
\includegraphics[height=15mm]{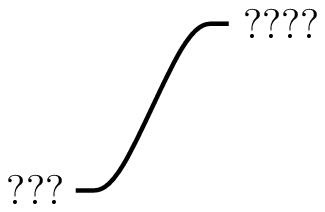}

\subsection*{Доказательство Гипотезы. Шаг 4. Пример}

Понятно, что процедуру {\it сокращения разворотов} можно применять, пока не останутся только развороты в одном направлении.
Попробуем сделать это с нашей кривой $\alpha'$:

\ver\hfil
\includegraphics[height=23mm]{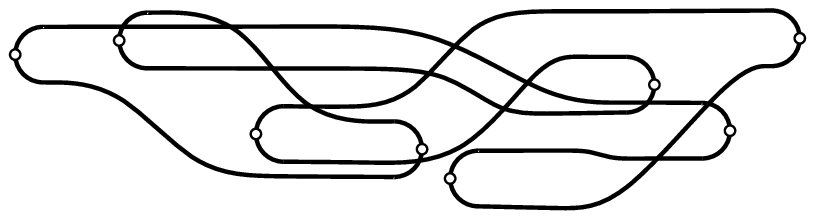}

\ver\hfil
\includegraphics[height=23mm]{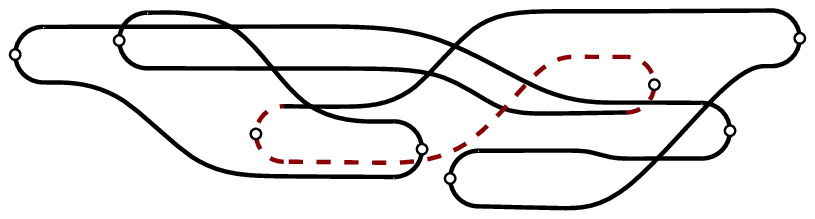}

\ver\hfil
\includegraphics[height=23mm]{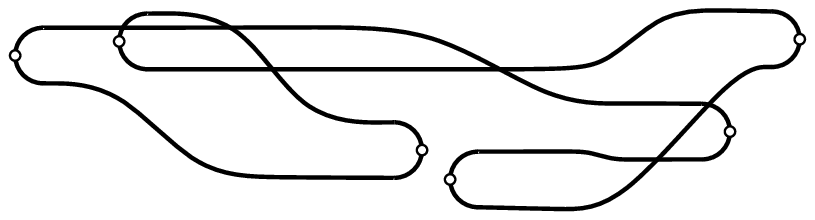}

\ver\hfil
\includegraphics[height=23mm]{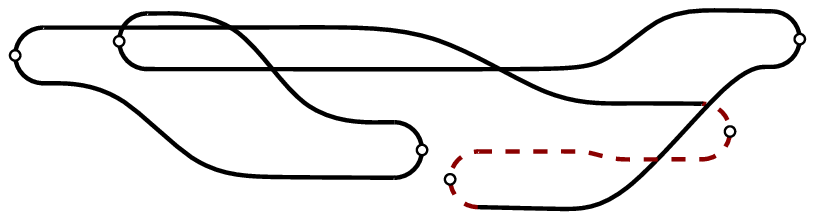}

\ver\hfil
\includegraphics[height=23mm]{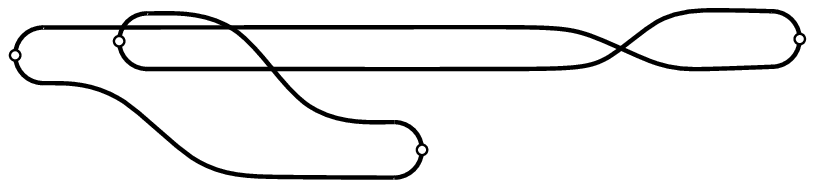}

\ver\hfil
\includegraphics[height=23mm]{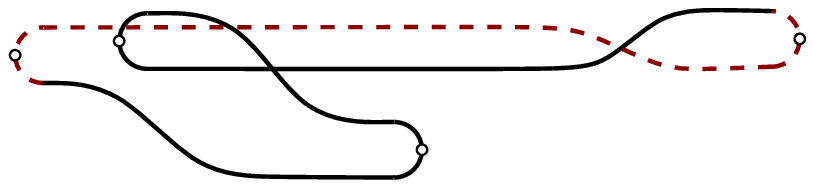}

\ver\hfil
\includegraphics[height=23mm]{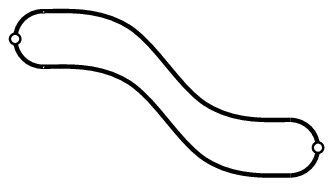}

\subsection*{Доказательство Гипотезы. Шаг 4}

Итак, процедуру {\it сокращения разворотов} можно применять, пока не останутся только развороты в одном направлении.
В итоге кривая будет выглядеть как-то так:

\ver\hfil
\includegraphics[height=35mm]{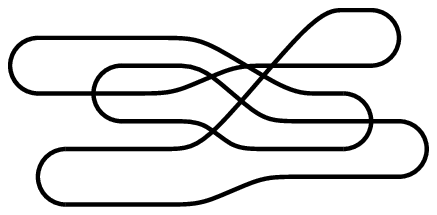}
\ver

Её можно аккуратно ``собрать'' в соответствующую кривую из таблицы.

\ver\hfil
\includegraphics[height=25mm]{curve-uncomplicated-1.eps}
\hfil
\includegraphics[height=25mm]{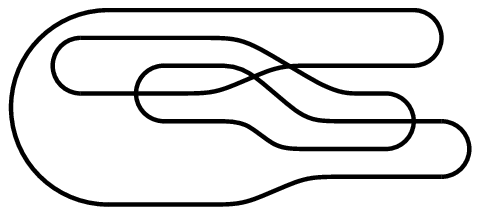}
\hfil
\includegraphics[height=25mm]{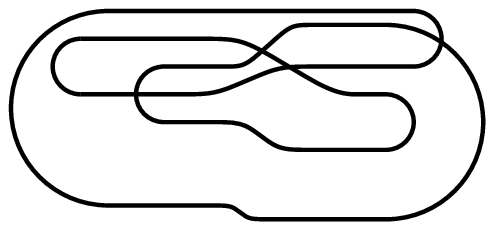}

\ver\hfil
\includegraphics[height=25mm]{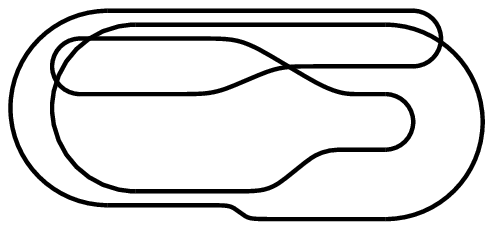}
\hfil
\includegraphics[height=25mm]{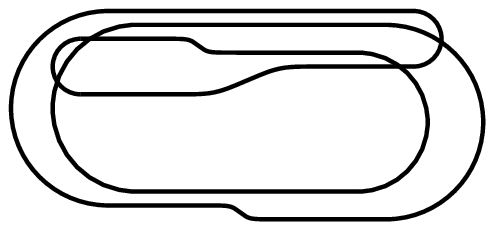}
\hfil
\includegraphics[height=25mm]{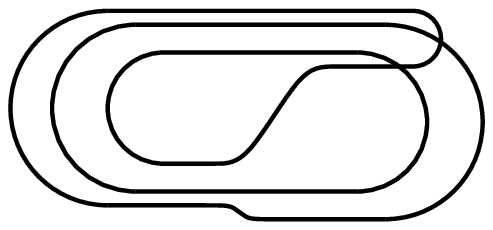}

\ver\hfil
\includegraphics[height=25mm]{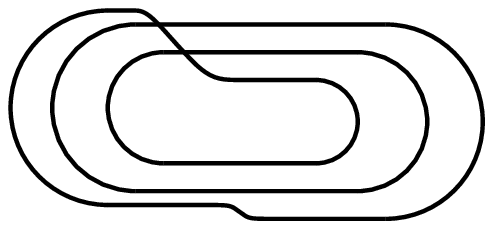}
\hfil
\includegraphics[height=25mm]{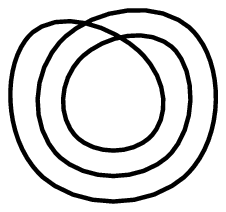}

\subsection*{Доказательство Гипотезы. Шаг 4. Финиш}

Итак, {\bf гипотеза доказана:}
мы умеем деформировать любую кривую $\alpha$ в кривую $\gamma_k$ для какого-нибудь~$k$.

Но погодите!
Описанная процедура не может дать в конце кривую $\gamma_0$!
Действительно, у полученной после сокращения кривой {\it все развороты идут в одном и том же направлении}, а у $\gamma_0$ --- нет.

\ver\hfil
\includegraphics[height=15mm]{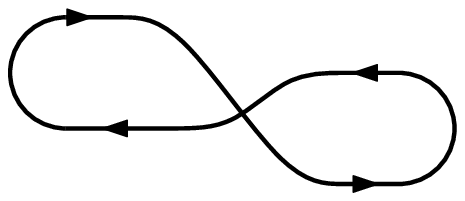}

Другими словами, даже если бы $\alpha$ была регулярно гомотопна $\gamma_0$,
в конце нашего алгоритма мы получаем $\gamma_k$ при некотором $k\ne0$.
Но мы уже доказали, что $\gamma_0$ не регулярно гомотопна $\gamma_k$,
что даёт противоречие.

\subsection*{Доказательство Гипотезы. Шаг 4. Уточнение}

На самом деле, на Шаге 3 процедуру ``потягивания'' за концы можно применять не ко всем кривым.
Действительно, в уточнении к Шагу 3 мы неявно используем, что фрагменты за потягиваемыми концами не соединяются.

Но если бы они соединялись, потянуть в разные стороны за концы было бы невозможно.

\ver\hfil
\includegraphics[height=18mm]{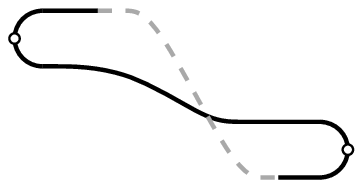}
\ver

\noindent
А это имеет место ровно в том случае, когда наша кривая --- это $\gamma_0$, либо $\gamma_0'$.

\subsection*{Доказательство Гипотезы. Шаг 4. Теперь точно финиш}

Таким образом, алгоритм завершается,
если

$\bullet$ если остались только развороты в одинаковом направлении;

$\bullet$ если остались два разворота в противоположных направлениях.

\noindent
Во втором случае получилась кривая $\gamma_0$ или $\gamma_0'$,
а в первом --- кривая $\gamma_k$ для $k\ne0$.

\subsection*{Доказательство Гипотезы. Замечание}

Вспомним, для какого $k\in\Z$ кривая $\alpha'$ регулярно гомотопна $\gamma_k$?
Разумеется, только для $k=r(\alpha)$.
Оказывается, конце Шага 2 можно вычислить значение $k$ наглядно.
Как это сделать?

\ver\hfil
\includegraphics[height=18mm]{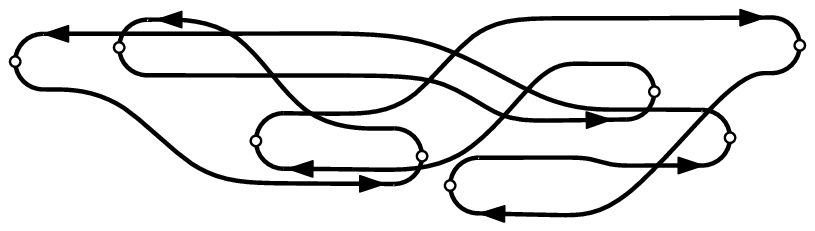}

\begin{problem}
Это число $k$ равно полуразности числа поворотов направо и числа поворотов налево,
то есть $$r(\alpha')=\dfrac{r_L-r_R}2.$$
\end{problem}

Например, для кривой на рисунке выше оно равно $\frac{5-3}2=1$.
И действительно, явный подсчёт показывает, что
при обходе кривой в заданном направлении
вектор скорости совершает суммарно 1 оборот против часовой стрелки.
Это полностью согласуется с разобранным выше примером к Шагу 4,
где мы показали что данная кривая регулярно гомотопна кривой $\gamma_1$.

\subsection*{Доказательство Гипотезы. Подведение итогов}

Можно резюмировать результат следующим образом.

\begin{theorem}[Уитни-Грауштейн]
Любая гладкая замкнутая кривая $\alpha$ на плоскости регулярно гомотопна ровно одной кривой
$\gamma_k$ из набора в \S\ref{s:collection}.
При этом номер $k$ равен числу оборотов вектора скорости~$r(\alpha)$.
\end{theorem}

Представленное доказательство является авторским,
проделанные в нём комбинаторно-геоме\-три\-чес\-кие представления кажутся довольно естественными.%
\footnote{Отдельный интерес представляет то, что в нашем доказательстве явно видны места,
где для аккуратного обоснования построений применимы классические теоремы анализа --- лемма Сарда.}
Два других доказательства Теоремы Уитни-Грау\-штей\-на можно найти в \S5.2 книги Прасолова~\cite{prasolov-ekd}.

\section{Обобщение Теоремы Уитни-Грауштейна}

В этом разделе мы неформально и в общих чертах наметим аналог доказанной выше теоремы для произвольных размерностей.
Обрисованные здесь определения не претендуют на полноту --- на полноценное их освоение уходит не один семестровый курс\ldots

\subsection{Многообразия}

{\it Гладкое многообразие} --- это, грубо говоря, множество,
склеенное из кусков евклидова пространства как папье-маше.
Куски евклидова пространства называются {\it картами}.
При этом функции склейки (отображения замены координат) между картами
должны быть непрерывно-дифференцируемыми (обычно --- бесконечное число раз).

У многообразия $X$ в каждой точке можно откладывать {\it касательные вектора},
про которые принято думать как про инфинитезимально-малые вектора в локальных координатах.
Множество всех касательных векторов во всех точках образует {\it касательное расслоение},
оно обозначается $TX$.

Пусть $M,N$ --- пара многообразий.
Для отображения $M\to N$ можно говорить про {\it производную}.
Более точно, если отображение $f:M\to N$ {\it гладкое} (непрерывно-дифференцируемое в локальных координатах),
для него определён {\it дифференциал} $df:TM\to TN$.

По сути, дифференциал отображения --- это прямое обобщение производной функции.
Дифференциал показывает, куда перешёл какой касательный вектор.
Наглядно: если точка $x$ движется в $M$ со скоростью $v$,
то по определению её образ $f(x)$ будет двигаться в $N$ со скоростью $df(v)$
в тот же момент времени.

Важное свойство дифференциала: для каждого слоя касательного расслоения к $M$,
т.\,е.\ множества векторов, выходящих из одной точки $M$,
дифференциал является {\it линейным отображением} (сохраняет структуру векторного пространства).

\subsection{Погружения}

Пусть $\dim M<\dim N$.
Гладкое отображение $M\to N$ называется {\it погружением},
если в каждой точке $x\in M$ дифференциал $df_x$ имеет ранг $\dim M$.

\begin{theorem}[Смейла-Хирша.]
Непрерывное отображение $f:M\to N$
можно продеформировать в погружение 
если и только если
существует послойное вложение $F:TM\to TN$, накрывающее~$f$.
\end{theorem}

Например, эта теорема позволяет легко доказать, что
существует погружение проективной плоскости в пространство.
Одно из таких погружений называется {\it поверхность Боя  (Boy surface)},
оно показано на рисунке ниже.

\ver\hfil
\includegraphics[height=60mm]{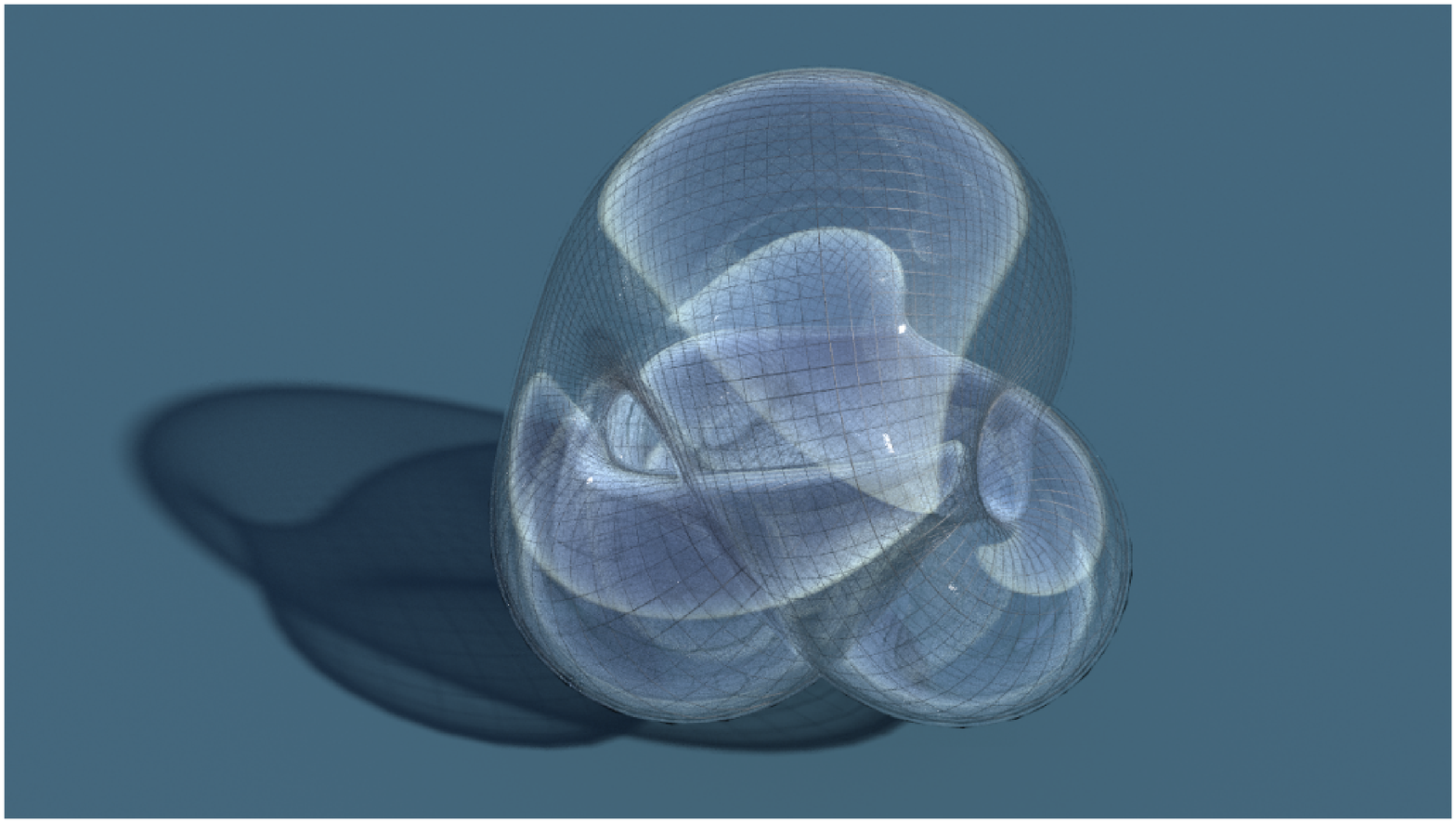}
\ver

Теорема Смейла-Хирша была доказана примерно в 60-х годах прошлого века и неоднократно усовершенствовалась.
В дальнейшем Миша Громов развил аналогичную технику для большего класса задач.
Эта техника называется {\it $h$-принцип} и активно применяется в современной науке.

\section{Вместо послесловия: упражнения для знатоков}

\begin{problem}
Докажите, что любую гладкую кривую можно разбить на конечное количество несамопересекающихся дуг.
\end{problem}

\begin{problem}
Верно ли, что любая гладкая кривая делит плоскость на конечное число частей?
\end{problem}

\begin{problem}
В определении из \S\ref{s:deform} опустим требование {\it непрерывности} изменения вектора скорости в процессе регулярной гомотопии.
Как изменится описание классов регулярной гомотопности замкнутых кривых на плоскости?
\end{problem}

\begin{problem}
Докажите, что для любой прямой $l$ и гладкой замкнутой кривой $\alpha$ существует по меньшей мере две точки на $\alpha$, касательные в которых перпендикулярны $l$.
\end{problem}

\begin{problem}
Постройте послойное вложение $T\RP^2\to T\R^3$.
\end{problem}

\section*{Благодарности}

Я горячо благодарен организаторам фестиваля \cite{dgf}
(за их работу, за возможность выступить, и вообще желаю им процветания и успехов в их труде),
а также участникам фестиваля, без которых также не возникла бы вся эта деятельность.
Кроме того, я признателен Дарье Рыченковой за многочисленные стилистические замечания,
которые несомненно помогли сделать статью лучше,
а Сашечке Пилипюк --- за ряд полезных уточнений по содержанию текста
(все оставшиеся огрехи, как языковые, так и математические --- остаются полностью на совести автора).


\vspace{3em}

{\sl
\noindent
Андрей Рябичев
\nopagebreak

\noindent
НМУ, школа №\,179 г.\ Москвы
\nopagebreak

\noindent
E-mail: \texttt{ryabichev@179.ru}
}

\end{document}